\documentclass[12pt]{article}
\usepackage{setspace}
\usepackage{cite}

\usepackage{mathrsfs}
\usepackage{amsmath}
\usepackage{amssymb}
\usepackage{amsfonts}
\usepackage{amsmath,amsfonts,amssymb}
\usepackage{graphicx}
\usepackage{subfigure}
\usepackage{caption2}
\usepackage[rightcaption]{sidecap}
\usepackage{fancyhdr}
\newcommand{\R}{\mathbb R}

\newcommand{\B}{\mathbb B}

\newtheorem{remark}{Remark}

\def \qed{\hfill{$\Box$}}

\begin{document}
\title{Singularity of the extremal solution for supercritical biharmonic equations with power-type
nonlinearity
}
\date{}
\author{Baishun Lai$^{1,2}$,  Zhengxiang Yan,  Yinghui Zhang\\
\small {\it 1 Institute of Contemporary Mathematics, Henan University}\\
\small {\it Kaifeng 475004, P.R.China}\\
\small {\it 2  School of Mathematics and Information Science,Henan University}\\
\small {\it  laibaishun@henu.edu.cn}\\
\small {\it Department of Mathematics, Hunan Institute of Science and Technology}\\
\small {\it   Yueyang, Hunan 414006, China}\\
 \small {\it Xinyang Vocational and Technical College, Xinyang 464000, China
}} \maketitle
\begin{center}
\begin{minipage}{120mm}
{\small { Abstract.}\ \ \ Let $\lambda^{*}>0$ denote the largest
possible value of $\lambda$ such that
$$
\left\{
\begin{array}{lllllll}
\Delta^{2}u=\lambda(1+u)^{p} &  \mbox{in}\ \ \B, \\
u=\frac{\partial u}{\partial n} =0 &  \mbox{on}\  \  \partial \B\\
\end{array}
\right.
$$
has a solution, where $\B$ is the unit ball in $R^{n}$ centered at the origin,
$p>\frac{n+4}{n-4}$ and $n$ is the exterior unit normal vector. We show that
for $\lambda=\lambda^{*}$ this problem possesses a unique weak solution $u^{*}$,
 called the extremal solution. We prove that $u^{*}$ is  singular when $n\geq 13$ for $p$ large enough, in which case
 $u^{*}(x)\leq r^{-\frac{4}{p-1}}-1$ on the unit ball and actually solve part of the open problem which
\cite{D} left.\vskip 0.1in


  }
\end{minipage}
\end{center}
\medskip
\noindent
 {\bf 1. Introduction and results} \medskip

\noindent

 In the previous two decades, positive solutions to the second order semilinear elliptic problem
$$
\left\{
\begin{array}{lllllll}
-\Delta u=\lambda g(u) &\ \ \mbox{in}\  \Omega,\\
u=0& \ \ \mbox{on}\ \Omega,
\end{array}
\right.
\eqno(1.1)
$$
have attracted a lot of interest, see e.g. \cite{Cr,Jo, Lion, Ma,Mi}
and references therein. Here, we only mention the work by Joseph and
Lundgren \cite{Jo}. In their well known work, Joseph and Lundgren
gave a complete characterization of all positive solutions of (1.1)
in the case $g(u)=e^{u}$ or $g(u)=(1+au)^{p}, ap>0$, $\lambda>0$ and
$\Omega$ is unit ball in $R^{n}$. In particular, they found a
remarkable phenomenon for $g(u)=e^{u}$ and $n>2$: either (1.1) has
at most one solution for each $\lambda$ or there is a value of
$\lambda$ for which infinitely many solutions exist. In the case of
a power nonlinearity the same alternative is valid if $n\geq3$ and
$p\not\in (1,\frac{n+2}{n-2}]$. As a subsequent step, P.L. Lions
(\cite{Lion}, section 4.2 (c)) suggests to study positive solutions
to systems of semilinear elliptic equations.
 So it is an important task to gain a deeper understanding for related higher order problems.\medskip

In this paper we study a semilinear equation involving the bilaplacian operator and a power type nonlinearity
$$
\left\{
\begin{array}{lllllll}
\Delta^{2}u=\lambda(1+u)^{p} &  \mbox{in}\ \ \B, \\
u=\frac{\partial u}{\partial n} =0 &  \mbox{on}\  \  \partial \B,\\
\end{array}
\right.
\eqno(1.2)
$$
where $\B\subset R^{n}$ is the unit ball, $\lambda>0$ is an
eigenvalue parameter, $n\geq5$ and $p\geq\frac{n+4}{n-4}$. The
subcritical case $p<\frac{n+4}{n-4}$ is by now \textquotedblleft
folklore",
 where existence and multiplicity results are easily established by means of variational methods. For the critical case $p=\frac{n+4}{n-4}$
  (under Navier boundary conditions), we refer to \cite{Ber}.
 Recently, a lot of research on supercritical case , i.e., $p>\frac{n+4}{n-4}$, has been done and many beautiful important results have been proved.
 In what follows, we will summarize some of the
results obtained by \cite{Fe,Fer}. For convenience, we introduce the
following notions:

\medskip \noindent \textit{{\bf Definition 1.1.} We
say that $u\in L^{p}(B)$ is a solution of (1.2) if $u\geq0$ and if
for all $\varphi\in C^{4}(\bar{\B})$ with $\varphi|\partial
\B=|\nabla \varphi||\partial \B=0$ one has
$$
\int_{\B}u\Delta^{2}\varphi dx=\lambda \int_{\B}(1+u)^{p}\varphi dx.
$$
We call u singular if $u \not\in L^{\infty}(\B)$, and regular if
$u\in L^{\infty}(\B)$.}

A radial singular solution $u=u(r)$ of (1.2) is called weakly
singular if $\lim_{r\to0}r^{\frac{4}{p-1}}u(r)\in [0, \infty]$
exists. \medskip

Note that by standard regularity theory for the biharmonic operator,
any regular solution $u$ of (1.2) satisfies $u\in
C^{\infty}(\bar{\B})$. Note also that by the positivity preserving
property of $\Delta^{2}$ in the ball (see \cite{Bo}) any solution of
(1.2) is positive, see also \cite{Ar1} for a generalized statement.
This property is known to fail in general domains. For this reason,
we restrict ourselves to ball. Hence, the sub-and super-solution
method applies as well as monotone iterative procedures.

\medskip
\noindent \textit{{\bf Definition 1.2.}  We call a solution $u$ of
(1.2) minimal if $u\leq v$ a.e. in $\B$ for any further solution $v$
of (1.2)} \medskip

We also denote by $\lambda_{1}>0$ the first eigenvalue for the biharmonic operator with Dirichlet boundary conditions
$$
\left\{
\begin{array}{lllllll}
\Delta^{2} \varphi_{1}=\lambda_{1}\varphi_{1} & \mbox{in}\ \ \B,\\
\varphi_{1}=\frac{\partial\varphi_{1}}{\partial n}=0& \mbox{on}\ \ \partial\B.
\end{array}
\right.
$$
It is known from the positivity preserving property and Jentzsch's
(or Krein-Rutman's) theorem that $\lambda_{1}$ is isolated and
simple and the corresponding eigenfunction $\varphi_{1}$ don't
change sign. \medskip

\noindent \textit{{\bf Definition 1.3.}\ We say a weak solution of
(1.2) is stable (resp. semi-stable) if
$$
\mu_{1}(u)=\inf\{\int_{\B}(\Delta \varphi)^{2}-p\lambda\int_{\B}\varphi^{2}(1+u)^{p-1}:
\phi\in H_{0}^{2}(\B), \|\phi\|_{L^{2}}=1\}
$$
is positive (resp. non-negative).}\medskip

To illuminate the motivations of this paper in detail, we need the following notations which
will be used throughout the paper. Set
$$
  K_{0}=\frac{8(p+1)}{p-1}\left[n-\frac{2(p+1)}{p-1}\right]\left[n-\frac{4p}{p-1}\right],
 $$
and
$$
p_{c}=\frac{n+2-\sqrt{4+n^{2}-4\sqrt{n^{2}+H_{n}}}}{n-6-\sqrt{4+n^{2}-4\sqrt{n^{2}+H_{n}}}} \ \ \quad \mbox{for}\ \  n\geq 3,
$$
with $H_{n}=(n(n-4)/4)^{2}$ and the number $p_{c}$ such that when $p=p_{c}$  then
$$
(\frac{4}{p-1}+4)(\frac{4}{p-1}+2)(n-2-\frac{4}{p-1})(n-4-\frac{4}{p-1})=H_{n}.
$$

Now we summarize some of the well-known results as follows:

\medskip
\noindent \textit{{\bf Theorem A}\cite{Fe, Fer}. There exists
$\lambda^{*} \in \left[K_{0}, \frac{\lambda_{1}}{p}\right)$ such
that: \vskip0.1in \noindent(i) For $\lambda\in(0,\lambda^{*})$,
(1.2) admits a minimal stable regular solution, denoted by
$u_{\lambda}$. This solution is radially symmetric and strictly
decreasing in $r=|x|$. \vskip0.1in \noindent (ii) For
$\lambda=\lambda^{*}$, (1.2) admits at least one not necessarily
bounded solution, which is called extremal solution $u^{*}$
\vskip0.1in\noindent(iii) For $\lambda>\lambda^{*}$, (1.2) admits no
(not even singular ) solutions } \medskip

\noindent \textit{{\bf Theorem B}\cite{D}.  Assume that
   $$
   \frac{n+4}{n-4}<p<p_{c}\ \  \mbox{if}\ \  n\geq 13, \ \  \frac{n+4}{n-4}<p<\infty \ \  \mbox{if}\ \ 5\leq n\leq12
   $$
Then, $u^{*}$ is regular.}\medskip

From Theorem B,  we know that the extremal solution of  $(1.2)$ is
regular for a certain range of $p$ and $n$.  At the same time, they
left a open  problem: if
$$
n\geq 13\ \ \mbox{and}\ \ p\geq p_{c},
$$
is $u^{*}$ singular?

In this paper, by constructing a semi-stable singular $H_{0}^{2}(\B)-$ weak sub-solution of (1.2), we prove that, if $p$ is large enough,
 the extremal solution is singular for dimensions $n\geq 13$ and {\bf complete part of the above open problem}. Our result is stated as follows:
\vskip 0.2in

\noindent \textit{{\bf  Theorem 1.1.} There exists $p_0 > 1$ large
enough such that for $p \geq p_0$, the unique extremal solution of
(1.2) is singular for dimensions $n\geq13$, in which case $u^{*}\leq
|x|^{-\frac{4}{p-1}}-1$ on the unit ball.} \vskip 0.2in

From the technical point of view,  one of the obstacle is the
well-known difficulty of extracting energy estimates for solutions
of fourth order problems from their stability properties. Besides,
for the corresponding second order problem (1.1), the starting point
was an explicit singular solution for a suitable eigenvalue
parameter $\lambda$ which turned out to  play a fundamental role for
the shape of the corresponding bifurcation diagram, see \cite{Br}.
When turning to the biharmonic problem (1.2) the second boundary
condition $\frac{\partial u}{\partial n}=0$ prevents to find an
explicit singular solution. This means that the method used to
analyze the regularity of the extremal solution for second order
problem could not carry to the corresponding problem for (1.2). In
this paper, we, in order to overcome the second obstacle, use
improved and non standard Hardy-Rellich inequalities recently
established by Ghoussoub-Moradifam in \cite{Gh} to construct a
semi-stable singular $H_{0}^{2}(\B)-$ weak sub-solution of (1.2).
\vskip 0.1in

This paper is organized as follows. In the next section, some preliminaries are reviewed.
 In Section 3, we will show that the extremal solution $u^{*}$ in dimensions $n\geq 13$ is singular by constructing a semi-stable singular
 $H_{0}^{2}(\B)-$ weak sub-solution of (1.2).
 \medskip

\noindent {\bf 2. Preliminaries}

\medskip First we give some comparison
principles which will be used throughout the paper.

\medskip
\noindent \textit{{\bf Lemma 2.1.}\  (Boggio's principle, \cite{Bo})
If $u\in C^{4}(\bar{\B}_{R})$ satisfies
$$
\left\{
\begin{array}{lllllll}
\Delta^{2}u\geq 0 & \mbox{in}\ \  \B_{R},\\
u=\frac{\partial u}{\partial n}=0 &  \mbox{on}\ \  \partial \B_{R},
\end{array}
\right.
$$
then $u\geq 0$ in $\B_{R}$.}

\medskip \noindent \textit{{\bf Lemma 2.2.}\
Let $u\in L^{1}(\B_{R})$ and suppose that
$$
\int_{\B_{R}}u\Delta^{2}\varphi\geq0
$$
for all $\varphi\in C^{4}(\bar{\B}_{R})$ such that $\varphi\geq0$ in $\B_{R}$,
$\varphi|_{\partial \B_{R}}=\frac{\partial \varphi}{\partial n}|_{\partial \B_{R}}=0$. Then $u\geq 0$ in $\B_{R}$.
 Moreover $u\equiv 0$ or $u>0$ a.e., in $\B_{R}$.} \medskip

\noindent For a proof see Lemma 17 in \cite{Ar1}.
\medskip

\noindent\textit{{\bf Lemma 2.3.}\ If $u\in H^{2}(\B_{R})$ is
radial, $\Delta^{2}u\geq 0$ in $\B_{R}$ in the weak sense, that is
$$
\int_{\B_{R}}\Delta u\Delta\varphi\geq 0 \ \ \forall \varphi \in C_{0}^{\infty}(\B_{R}), \ \varphi\geq0
$$
and $u|_{\partial \B_{R}}\geq0, \frac{\partial u}{\partial
n}|_{\partial \B_{R}}\leq 0$, then $u\geq 0$ in $\B_{R}$.}\medskip

\noindent{\bf Proof.} For the sake of completeness, we include a
brief proof here.  We only deal with the case $R=1$ for simplicity.
Solve
$$
\left\{
\begin{array}{lllllll}
\Delta^{2}u_{1}=\Delta^{2}u & \mbox{in} \ \  \B\\
u_{1}=\frac{\partial u_{1}}{\partial n}=0 &  \mbox{on}\ \ \partial \B
\end{array}
\right.
$$
in the sense $u_{1}\in H_{0}^{2}(\B)$ and $\int_{\B}\Delta u_{1}\Delta\varphi=\int_{\B}\Delta u\Delta\varphi$ for all $\varphi \in C_{0}^{\infty}(\B)$.
 Then $u_{1}\geq 0$ in $\B$ by lemma 2.2.

Let $u_{2}=u-u_{1}$ so that $\Delta^{2}u_{2}=0$ in $\B$. Define $f=\Delta u_{2}$. Then $\Delta f=0$ in $\B$ and since $f$ is radial we find that $f$
is a constant.
 It follows that $u_{2}=ar^{2}+b$. Using the boundary conditions we deduce $a+b\geq 0$ and $a\leq0$, which imply $u_{2}\geq0$.\vskip 0.1in

Now we give a notion of $H_{0}^{2}$($\B$)- weak solutions, which is
an intermediate class between classical and weak solutions.

\medskip \noindent
 \textit{{\bf Definition 2.1.}\ \  We say that $u$ is a
$H_{0}^{2}$($\B$)- weak solution of (1.2) if  $(1+u)^{p}\in
L^{1}(\B)$ and if
$$
\int_{\B}\Delta u\Delta\phi=\lambda\int_{\B}\phi(1+u)^{p},\ \ \ \forall \phi\in C^{4}(\bar{\B})\cap H_{0}^{2}(\B).
$$
 We say that $u$ is a $H_{0}^{2}$($\B$)- weak super-solution (resp. $H_{0}^{2}$($\B$)- weak sub-solution) of (1.2) if for $\phi\geq0$
 the equality is replaced with $\geq$ (resp.$\leq$) and $u\geq 0$ (resp. $\leq$), $\frac{\partial u}{\partial n}\leq 0$ (resp. $\geq$)
 on $\partial \B$.}\medskip

\noindent We also need the following comparison principle.

\medskip \noindent
\textit{{\bf Lemma 2.4.}\  \  Let  $u_{1}, u_{2}\in H^{2}(\B_{R})$
with $(1+u_{1})^{p}, (1+u_{2})^{p}\in L^{1}(\B_{R})$. Assume that
$u_{1}$ is stable and
$$
\Delta^{2}u_{1}\leq\lambda (1+u_{1})^{p}\ \ \mbox{in}\ \ \B_{R}
$$
in the $H^{2}(\B_{R})-$ weak sense, i.e.,
$$
\int_{\B_{R}}\Delta u_{1}\Delta\phi\leq\lambda\int_{\B_{R}}(1+u_{1})^{p}\phi, \ \ \forall \phi\in C_{0}^{\infty}(\B_{R}), \phi\geq0 \eqno(2.1)
$$
and $\Delta^{2}u_{2}\geq\lambda (1+u_{2})^{p}\ \ \mbox{in}\  \B_{R}$ in the similar weak sense. Suppose also
$$
u_{1}|_{\partial\B_{R}}=u_{2}|_{\partial\B_{R}}\ \ \mbox{and}\ \ \frac{\partial u_{1}}{\partial n}|_{\partial \B_{R}}=
\frac{\partial u_{2}}{\partial n}|_{\partial \B_{R}}.
$$}

\textit{Then
$$
 u_{1}\leq u_{2}  \ \ \mbox{in}\ \  \B_{R}.
$$}

\vskip 0.1in \noindent{\bf Proof.}  Define $\omega:=u_{1}-u_{2}$.
Then by the Moreau decomposition \cite{Mo} for the biharmonic
operator, there exist $\omega_{1}, \omega_{2}\in H_{0}^{2}(\B_{R})$,
with $\omega=\omega_{1}+\omega_{2}, \omega_{1}\geq0$ a.e.,
$\Delta^{2}\omega_{2}\leq0$ in the $H_{0}^{2}(\B_{R})-$ weak sense
and
$$
\int_{\B_{R}}\Delta \omega_{1}\Delta \omega_{2}=0.
$$
By Lemma 1.1, we have that $\omega_{2}\leq 0$ a.e. in $\B_{R}$.

Given now $0\leq\varphi\in C_{0}^{\infty}(\B_{R})$, we have that
$$
\int_{\B_{R}}\Delta\omega\Delta\varphi\leq\lambda\int_{\B_{R}}(f(u_{1})-f(u_{2}))\varphi,
$$
where $f(u)=(1+u)^{p}$. Since $u$ is semi-stable and by density one
has
$$
\lambda\int_{\B_{R}}f'(u)\omega_{1}^{2}\leq \lambda\int_{\B_{R}}(\Delta\omega_{1})^{2}=\lambda\int_{\B_{R}}\Delta\omega\Delta\omega_{1}
\leq\lambda\int_{\B_{R}}(f(u_{1})-f(u_{2}))\omega_{1}.
$$
Since $\omega_{1}\geq \omega$, one also has
$$
\int_{\B_{R}}f'(u)\omega \omega_{1}\leq \int_{\B_{R}}(f(u_{1})-f(u_{2}))\omega_{1}
$$
which once re-arrange gives
$$
\int_{\B_{R}}\tilde{f}\omega_{1}\geq 0,
$$
where $\tilde{f}(u_{1})=f(u_{1})-f(u_{2})-f'(u_{1})(u_{1}-u_{2})$.
The strict convexity of $f$ gives $\tilde{f}\leq0$ and $\tilde{f}<0$
whenever $u\neq U$. Since $\omega_{1}\geq0$ a.e. in $\B_{R}$, one
sees that $\omega\leq0$ a.e. in $\B_{R}$. The inequality $u_{1}\leq
u_{2}$ a.e. in $\B_{R}$ is then established.\medskip

\noindent The following variant of lemma 2.4 also holds:

\medskip\noindent \textit{{\bf Lemma 2.5.}\ \ Let  $u_{1}, u_{2}\in H^{2}(\B_{R})$
be radial with $(1+u_{1})^{p}, (1+u_{2})^{p}\in L^{1}(\B_{R})$.
Assume $ \Delta^{2}u_{1}\leq\lambda (1+u_{1})^{p}\ \ \mbox{in}\
\B_{R}$ in the sense of (2.1) and and $\Delta^{2}u_{2}\geq\lambda
(1+u_{2})^{p}\ \ \mbox{in}\  \B_{R}$. Suppose $u_{1}|_{\partial
\B_{R}}\leq u_{2}|_{\partial \B_{R}}$ and $\frac{\partial
u_{1}}{\partial n}|_{\partial \B_{R}}\geq \frac{\partial
u_{2}}{\partial n}|_{\partial \B_{R}}$ and suppose also that $u_{1}$
is semi-stable. Then $u_{1}\leq u_{2}$ in $\B_{R}$.}\medskip

\noindent{\bf Proof.} We solve for $\hat{u}\in H_{0}^{2}(\B)$ such
that
$$
\int_{\B_{R}}\Delta \hat{u}\Delta\phi=\int_{\B_{R}}\Delta(u_{1}-u_{2})\Delta\phi\ \ \forall \phi\in C_{0}^{\infty}(\B_{R}).
$$
By Lemma 2.3 it follows that $\hat{u}\geq u_{1}-u_{2}$. Next we apply the Moreau decomposition to $\hat{u}$, that is
$\hat{u}=w+v$ with $w,v \in H_{0}^{2}(\B_{R}), w\geq0, \Delta^{2}v\leq0$ in $\B_{R}$ and $\int_{\B_{R}}\Delta w\Delta v=0$.
 Then the argument follows that of Lemma 2.4.

\vskip 0.1in \noindent \textit{{\bf Lemma 2.6.}\ \ Let $u$ be a
semi-stable $H^{2}_{0}(\B)-$ weak solution of (1.2). Assume $U$ is a
$H^{2}_{0}(\B)-$ super-solution of (1.2). Then if $u$ is a classical
solution and $\mu_{1}(u)=0$, we have $u=U$.}

\medskip
 \noindent{\bf Proof.} Since $u$ is a classical solution, it is
easy to see that the infimum in $\mu_{1}(u)$ is attained at some
$\varphi$. The function $\varphi$ is then the first eigenfunction of
$\Delta^{2}-\lambda f'(u)$ in $H_{0}^{2}(\B)$, where
$f(u)=(1+u)^{p}$. Now we show that $\phi$ is of fixed sign.
 Using the Moreau decomposition, one has $\phi=\phi_{1}+\phi_{2}$ where $\phi_{i}\in H_{0}^{2}(\B)$ for $i=1,2, \phi_{1}\geq0,
\int_{\B}\Delta\phi_{1}\Delta\phi_{2}=0$ and $\Delta^{2}\phi_{2}\leq0$ in the $H_{0}^{2}(\B)-$ weak sense. If $\phi$ changes sign, then
$\phi_{1}\not\equiv0$ and $\phi_{2}<0$ in $\B$. We can write now:
$$
0=\mu_{1}(u)\leq\frac{\int_{\B}(\Delta(\phi_{1}-\phi_{2}))^{2}-\lambda f'(u)(\phi_{1}-\phi_{2})^{2}}{\int_{\B}(\phi_{1}-\phi_{2})^{2}}
<\frac{\int_{\B}(\Delta \phi)^{2}-\lambda f'(u)\phi^{2}}{\int_{\B}\phi^{2}}=\mu_{1}(u)
$$
in view of $\phi_{1}\phi_{2}<-\phi_{1}\phi_{2}$ in a set of positive measure, leading to a contradiction.

So we can assume $\phi\geq0$, and by the Boggio's principle we have $\phi>0$ in $\B$. For $0\leq t\leq 1$ define
$$
g(t)=\int_{\B}\Delta(tU+(1-t)u)\Delta\phi-\lambda\int_{\B}f(tU+(1-t)u)\phi,
$$
where $\phi$ is the above first eigenfunction. Since $f$ is convex one sees that
$$
g(t)\geq\lambda \int_{\B}[tf(U)+(1-t)f(u)-f(tU+(1-t)u)]\phi\geq0
$$
for every $t\geq0$. Since $g(0)=0$ and
$$
g'(0)=\int_{\B}\Delta(U-u)\Delta\phi-\lambda f'(u)(U-u)\phi=0
$$
we get that
$$
g''(0)=-\lambda\int_{\B}f''(u)(U-u)^{2}\phi\geq0.
$$
Since $f''(u)\phi> 0$ in $\B$, we finally get that $U=u$ a.e. in
$\B$. \medskip

 \noindent From this lemma, we immediately obtain: \medskip

\noindent \textit{{\bf Corollary 2.1} (i)\  When $u^{*}$ is a
classical solution, then $\mu_{1}(u^{*})=0$ and $u^{*}$ is the
unique $H^{2}_{0}(\B)-$ weak solution of (1.2); \vskip0.1in
\noindent(ii)\ If $v$ is a singular semi-stable $H^{2}_{0}(\B)-$
weak solution of (1.2), then $v=u^{*}$ and $\lambda=\lambda^{*}$.}

\medskip
\noindent {\bf Proof.} (i) Since the function $u^{*}$ is a classical
solution, and by the Implicit Function Theorem we have that
$\mu_{1}(u^{*})=0$ to prevent the continuation of the minimal branch
beyond $\lambda^{*}$. By Lemma 2.4, $u^{*}$ is then the unique
$H_{0}^{2}(\B)-$ weak solution of (1.2).

(ii) Assume now that $v$ is a singular semi-stable $H_{0}^{2}(\B)-$
weak solution of (1.2). If $\lambda<\lambda^{*}$, then by the
uniqueness of the semi-stable solution, we have $v=u_{\lambda}$. So
$v$ is not singular and a contradiction arises. By Theorem A (iii)
we have that $\lambda=\lambda^{*}$. Since $v$ is a semi-stable
$H^{2}_{0}(\B)-$ weak solution of (1.2) and $u^{*}$ is a
$H^{2}_{0}(\B)-$ weak super-solution of (1.2), we can apply Lemma
2.4 to get $v\leq u^{*}$ a.e. in $\B$. Since $u^{*}$ is a
semi-stable solution too, we can reverse the roles of $v$ and
$u^{*}$ in Lemma 2.5 to see that $v\geq u^{*}$ a.e. in $\B$. So
equality $v=u^{*}$ holds and the proof is complete. \qed

\medskip\noindent {\bf 3 Proof of Theorem 1.1}\medskip

\noindent Inspired by the work of \cite{Mo1},  we will first show
the following upper bound on $u^{*}$

\medskip\noindent \textit{{\bf Lemma 3.1.} If $n\geq 13$ and $p>p_{c}$, then
$u^{*}\leq |x|^{-\frac{4}{p-1}}-1$.}\medskip

\noindent {\bf Proof.}  Recall from Theorem A that
$K_{0}<\lambda^{*}$. We now claim that $u_{\lambda}\leq
\tilde{u}:=|x|^{-\frac{4}{p-1}}-1$ for all $\lambda\in(K_{0},
\lambda^{*})$. Indeed, fix such a $\lambda$ and assume by
contradiction that
$$
R_{1}:=\inf\{0\leq R\leq1: u_{\lambda}< \bar{u}\ \ \mbox{ in the interval}\ (R,1)\}>0.
$$
From the boundary conditions, one has that
$$
u_{\lambda}(r)<\tilde{u}(r) \ \ \mbox{as}\ \ r\to 1^{-}.
$$
Hence,
$$
0<R_{1}<1, u_{\lambda}(R_{1})=\tilde{u}(R_{1})\ \ \mbox{and}\ \ u_{\lambda}'(R_{1})\leq \tilde{u}'(R_{1}).
$$

Now consider the following problem
$$
\left\{
\begin{array}{lllllll}
\Delta^{2}u=K_{0}(1+u)^{p}& \ \mbox{in}\ \ \B_{R_{1}};\\
u=u_{\lambda}(R_{1})& \ \  \mbox{on}\ \ \partial \B_{R_{1}}; \\
 \frac{\partial u}{\partial n}= u'_{\lambda}(R_{1}) & \ \  \mbox{on}\ \ \partial \B_{R_{1}}.
\end{array}
\right.
$$
Then $u_{\lambda}$ is a super-solution to above problem while $\tilde{u}$ is a sub-solution to the same problem. Moreover for $n\geq 13$, we have
$$
pK_{0}\leq H_{n}:=\frac{n^{2}(n-4)^{2}}{16}
$$
and
$$
\int_{\B_{R_{1}}}(\Delta\phi)^{2}\geq  H_{n}\int_{\B_{R_{1}}}\frac{\phi^{2}}{|x|^{4}}dx\geq pK_{0}\int_{\B_{R_{1}}}(1+\tilde{u})^{p-1}.
$$
So $\tilde{u}$ is semi-stable and we deduce that $u_{\lambda}>\tilde{u}$ by the Lemma 2.4,  and a contradiction arises in view of the fact
$$
|u_{\lambda}|_{L^{\infty}(\B_{R_{1}})}<\infty,\ \ \mbox{and}\ \ |\tilde{u}|_{L^{\infty}(\B_{R_{1}})}=\infty.
$$
The proof is done.\qed
\medskip

In order to prove Theorem 1.1, we will need  a suitable
Hardy-Rellich type inequality which was established by
Ghoussoub-Moradifam in \cite{Gh}. It is stated as follows:

\medskip
\noindent \textit{{\bf Lemma 3.2.} Let $n\geq 5$ and $\B$ be the
unit ball in $\R^{n}$. Then there exists $C>0$, such that  the
following improved Hardy-Rellich inequality holds for all
$\varphi\in H_{0}^{2}(\B)$:
$$
\int_{\B}(\Delta \varphi)^{2}dx\geq \frac{n^{2}(n-4)^{2}}{16}\int_{\B}\frac{\varphi^{2}}{|x|^{4}}dx+C\int_{\B}\varphi^{2}dx.
$$}
\vskip 0.1in

\noindent \textit{{\bf Lemma 3.3.} Let $n\geq 5$ and $\B$ be the
unit ball in $\R^{n}$. Then the following improved Hardy-Rellich
inequality holds for all $\varphi\in H_{0}^{2}(\B)$:
\begin{eqnarray*}
\int_{\B}(\Delta \varphi)^{2} dx &\geq& \frac{(n-2)^{2}(n-4)^{2}}{16}\int_{\B}\frac{\varphi^{2}dx}{(|x|^{2}-0.9|x|^{\frac{n}{2}+1})(|x|^{2}-|x|^{\frac{n}{2}})}\\
&+& \frac{(n-1)(n-4)^{2}}{4}\int_{\B}\frac{\varphi^{2}dx}{|x|^{2}(|x|^{2}-|x|^{\frac{n}{2}})}. \hspace*{4.7cm} (3.0)
\end{eqnarray*}
As a consequence, the following improvement of the classical Hardy-Rellich inequality holds:
$$
\int_{\B}(\Delta \varphi)^{2} dx\geq \frac{n^{2}(n-4)^{2}}{16}\int_{\B}\frac{\varphi^{2}}{|x|^{2}(|x|^{2}-|x|^{\frac{n}{2}})}.
$$}

\medskip
\noindent We now give the following lemma which is crucial for the
proof of the Theorem 1.1. \medskip

\noindent {\bf Lemma 3.4.} Suppose there exist $\lambda'>0$ and a
radial function $u\in H^{2}(\B)\cap
W_{loc}^{4,\infty}(\B\setminus\{0\})$ such that $u\not\in
L^{\infty}(\B)$ and
\begin{align*}
\Delta^{2}u\leq\lambda'(1+u)^{p}\ \ \mbox{for}\  \ 0<r<1;\quad
u(1) = u'(1)=0
\end{align*}
and
\begin{align*}
p\beta\int_{\B}\varphi^{2}(1+u)^{p+1}\leq\int_{\B}(\Delta\varphi)^{2}\ \ \mbox{ for all}\  \varphi\in H_{0}^{2}(\B)
\end{align*}
for either $\beta>\lambda'$ or $\beta=\lambda'=\frac{H_{n}}{p}$. Then $u^{*}$ is singular and
$$
\lambda^{*}\leq \lambda'.
\eqno(3.1)
$$

\noindent {\bf Proof.} First, we prove $\lambda^{*}\leq \lambda'$.
Note that the stability and $u\in
L_{loc}^{\infty}(\B\setminus\{0\})$ yield to $(1+u)^{p}\in
L^{1}(\B)$, we easily see that $u$ is a weak sub-solution of
$(1.2)$. If now $\lambda'< \lambda^{*}$, by Lemma 2.5, $u$ would
necessarily be below the minimal solution $u_{\lambda'}$, which is a
contraction since $u$ is singular while $u_{\lambda'}$ is regular.

\medskip
Suppose first that $\beta=\lambda'=\frac{H_{n}}{p}$ and that $n\geq 13$. From the above we have $\lambda^{*}\leq\frac{H_{n}}{p}$,
we get from Lemma 3.1 and the improved Hardy-Rellich inequality that there exists $C>0$ so that for all $\phi\in H_{0}^{2}(\B)$
$$
\int_{\B}(\Delta \phi)^{2}-p\lambda^{*}\int_{\B}\phi^{2}(1+u^{*})^{p+1}\geq \int_{\B}(\Delta \phi)^{2}-
H_{n}\int_{\B}\frac{\phi^{2}}{|x|^{4}}\geq C\int_{\B}\phi^{2}
$$
It follows that $\mu_{1}(u^{*})>0$ and $u^{*}$ must therefore be singular since otherwise,
one could use the Implicit Function Theorem to continue the minimal branch beyond $\lambda^{*}$

\medskip
Suppose now that $\beta>\lambda'$ and let $\frac{\lambda'}{\beta_{1}}<\gamma<1$ and $\alpha:=(\frac{\gamma \lambda^{*}}{\lambda'})^{\frac{1}{p+1}}$
and define $\bar{u}:=\alpha^{-1}(1+u)-1$. We claim that
$$
u^{*}\leq \bar{u}\ \ \ \mbox{in}\  \B.
\eqno(3.2)
$$
To prove this, we shall show that for every $\lambda<\lambda^{*}$
$$
u_{\lambda}\leq \bar{u}\ \ \ \mbox{in}\  \B. \eqno(3.3)
$$
Indeed, we have
$$
\Delta^{2}\bar{u}=\alpha \Delta^{2}u\leq \alpha \lambda'(1+u)^{p}=\alpha^{p+1}\lambda'(1+\bar{u})^{p}.
$$
Now by the choice of $\alpha$, we have $\alpha^{p+1}\lambda'<\lambda^{*}$. To prove (3.3), it suffices to prove it
 for $\alpha^{p+1}\lambda'<\lambda<\lambda^{*}$.
Fix such $\lambda$ and assume that (3.3) is not true. Then
$$
\Lambda = \{ 0\leq R\leq 1| u_{\lambda}(R)>\bar{u}(R)\}
$$
is non-empty. Since $\bar{u}(1)=\alpha^{-1}-1>0=u_{\lambda}(1)$, we have $0<R_{1}<1, u_{\lambda}(R_{1})=\bar{u}(R_{1})$, and
$u'_{\lambda}(R_{1})\leq \bar{u}'(R_{1})$.
Now consider the following problem
$$
\left\{
\begin{array}{lllllll}
\Delta^{2}u=\lambda(1+u)^{p}& \ \mbox{in}\ \ \B_{R_{1}},\\
u=u_{\lambda}(R_{1})& \ \  \mbox{on}\ \ \partial \B_{R_{1}},\\
 \frac{\partial u}{\partial n}= u'_{\lambda}(R_{1}) & \ \  \mbox{on}\ \ \partial \B_{R_{1}}.
\end{array}
\right.
$$
Then $u_{\lambda}$ is a solution to above problem while $\bar{u}$ is a sub-solution to the same problem.
Moreover $\bar{u}$ is stable since $\lambda<\lambda^{*}$ and
$$
p\lambda(1+\bar{u})^{p+1}\leq p\lambda^{*}\alpha^{-(p+1)}(1+u)^{p+1}=p\lambda'\gamma^{-1}(1+u)^{p+1}
<p\beta_{1}(1+u)^{p+1},
$$
we deduce $\bar{u}\leq u_{\lambda}$ in $\B_{R_{1}}$, which is
impossible, since $\bar{u}$ is singular while $u_{\lambda}$ is
smooth. This establishes (3.2). From (3.2) and the above
inequalities, we have
$$
p\lambda^{*}(1+u^{*})^{p+1}\leq p\lambda'\gamma^{-1}(1+u)^{p+1}<p\beta_{1}(1+u)^{p+1}.
$$
Thus
$$
\inf_{\varphi\in C_{0}^{\infty}(\B)}\frac{\int_{\B}(\Delta \varphi)^{2}-p\lambda^{*}\varphi^{2}(1+u^{*})^{p+1}}{\int_{\B}\varphi^{2}}>0.
$$
This is not possible if $u^{*}$ is a smooth function by the Implicit Theorem. \qed

\medskip \noindent
 {\bf Proof Theorem 1.1}\ \  Uniqueness and the
upper bound estimate of the extremal solution $u^{*}$ have been
proven by Corollary 3.1 and Lemma 3.1.
   Now we only prove that $u^{*}$ is a singular solution of (1.1) for $n\geq 13$, in order to achieve this, we shall
find  a singular $H-$weak sub-solution of (1.1), denote by $\omega_{m}(r)$,  which is stable, according to the Lemma 3.4.\vskip 0.1in

Choosing
$$
\omega_{m}=a_{1}r^{-\frac{4}{p-1}}+a_{2}r^{m}-1, \ \ K_{0}=\frac{8(p+1)}{p-1}\left[n-\frac{2(p+1)}{p-1}\right]\left[n-\frac{4p}{p-1}\right]
$$
since $\omega(1)=\omega'(1)=0$, we have
$$
a_{1}=\frac{m}{m+\frac{4}{p-1}}\ \ \ \  a_{2}=\frac{\frac{4}{p-1}}{m+\frac{4}{p-1}}
$$
For any $m$ fixed, when $p\to+\infty$, we have
$$
a_{1}=1-\frac{4}{(p-1)m}+o(p^{-1}),\ \ \ a_{2}=1-a_{1}=\frac{4}{(p-1)m}+o(p^{-1})
$$
and
$$
K_{0}=\frac{8(n-2)(n-4)}{p}+o(p^{-1})
$$
Note that
\begin{eqnarray*}
&&\lambda'K_{0}(1+\omega_{m}(r))^{p}-\Delta^{2}\omega_{m}(r)=\lambda'K_{0}(1+\omega_{m}(r))^{p}
-a_{1}K_{0}r^{-\frac{4p}{p-1}}- a_{2}K_{1}r^{m-4}\\
&=&\lambda'K_{0}(a_{1}r^{-\frac{4}{p-1}}+a_{2}r^m)^{p}-a_{1}K_{0}r^{-\frac{4p}{p+1}}
-a_{2}K_{1}r^{m-4} \\
&=&K_{0}r^{-\frac{4p}{p-1}}\left[\lambda'(a_{1}+a_{2}r^{m+\frac{4}{p-1}})^{p}-a_{1}
-a_2K_{1}K_{0}^{-1}r^{\frac{4p}{p-1}+m-4}\right]\\
&=&K_{0}r^{-\frac{4p}{p-1}}\left[\lambda'(a_{1}+a_{2}r^{m+\frac{4}{p-1}})^{p}-a_{1}
-a_2K_{1}K_{0}^{-1}r^{m-\frac{4}{p-1}}\right]\\
&=&K_{0}r^{-\frac{4p}{p-1}}(a_{1}+a_{2}r^{m+\frac{4}{p-1}})^p\left[\lambda' - H(r^{m +\frac{4}{p-1}})\right] \hspace*{4.7cm} (3.4)\\
  \end{eqnarray*}
with
$$H(x) = (a_1+a_2x)^p\left[a_{1}+a_2K_{1}K_{0}^{-1}x\right], \  K_{1}=m(m-2)(m+n-2)(m+n-4) \eqno(3.5)$$

(1) Let $m=2$ and $n\geq 32$, then we can prove that
$$
\sup_{[0, 1]} H(x) = H(0) = a_1^{1-p}\longrightarrow e^{2}\ \  \mbox{as}\ \  p \longrightarrow +\infty.
 $$
So $(3.4)\geq0$ is valid as long as
$$\lambda' = e^{2}.$$
At the same time, we have ( since
$a_{1}+a_{2}r^{m+\frac{4}{p-1}}\leq a_{1}+a_{2}\leq 1$ in $[0,1]$)
$$
\frac{n^2(n-4)^{2}}{16}\frac{1}{r^{4}} -p\beta_{n}r^{-4}(a_{1}+a_{2}r^{2+\frac{4}{p-1}})^{p-1}
 \geq r^{-4}\left[\frac{n^2(n-4)^{2}}{16} - p\beta\right].\eqno(3.6)
$$
Let $\beta=(\lambda'+\varepsilon)K_{0}$, where $\varepsilon$ is arbitrary sufficient small,  we need finally here
$$\frac{n^2(n-4)^{2}}{16} - p\beta = \frac{n^2(n-4)^{2}}{16} - p(\lambda'+\varepsilon)K_{0} >0.$$
For that, it is sufficient to have for $p \longrightarrow +\infty$
$$\frac{n^2(n-4)^{2}}{16} - 8(e^{2}+\varepsilon)(n-2)(n-4)+o(\frac{1}{p})>0.$$
So $(3.6)\geq0$ holds only for $n\geq 32$ when $p\longrightarrow +\infty$. Moreover, for $p$ large enough
$$
8e^{2}(n-2)(n-4)\int_{\B}\varphi^{2}(1+\omega_{2})^{p+1}\leq
H_{n}\int_{\B}\frac{\varphi^{2}}{|x|^{4}}\leq \int_{\B}|\Delta \varphi|^{2}
$$
Thus it follows from Lemma 3.4 that $u^{*}$ is singular with $\lambda'=e^{2}K_{0}, \beta=(e^{2}K_{0}+\varepsilon(n,p))$
and $\lambda^{*}\leq e^{2}K_{0}$.

(2) Assume $13\leq n\leq 31$. We shall show that $u=\omega_{3.5}$ satisfies the assumptions of Lemma 5.4 for each dimension $13\leq n\leq 31$.
Using Maple, for each dimension $13\leq n\leq31$ one can verify that inequality $(3.4)\geq0$ holds for the $\lambda'$ given by Table 1.
Then, by using Maple again, we show that there exists $\beta>\lambda'$ such that
\begin{eqnarray*}
\frac{(n-2)^2(n-4)^2}{16}\frac{1}{(|x|^{2}-0.9|x|^{\frac{n}{2}+1})(|x|^{2}-|x|^{\frac{n}{2}})}\\
+\frac{(n-1)(n-4)^{2}}{4}\frac{1}{|x|^{2}(|x|^{2}-|x|^{\frac{n}{2}})}&\geq& p\beta(1+w_{3.5})^{p+1}.
\end{eqnarray*}
The above inequality and and improved Hardy-Rellich inequality (5.0)
guarantee that the stability condition (5.2) holds for
$\beta>\lambda'$. Hence by Lemma 3.4 the extremal solution is
singular for $13\leq n\leq31$ the value of $\lambda'$ and $\beta$
are shown in Table 1.\medskip

\begin{remark}
The values of $\lambda'$ and $\beta$ in Table 1 are not optimal.
\end{remark}

\begin{remark}
The improved Hardy-Rellich inequality (3.0) is crucial to prove that $u^{*}$ is singular in dimensions $n\geq 13$. Indeed
by the classical Hardy-Rellich inequality and $u:=w_{2}$, Lemma 5.4 only implies that $u^{*}$ is singular n dimensions $n\geq 32$.
\end{remark}

\noindent
 \textbf{Acknowledgements.}
The first author  would like to thank his advisor Prof. Yi-Li
 for his constant support, and encouragement. This research is supported in part  by National Natural Science
Foundation of China (Grant No. 10971061).


\newpage

$$
\begin{tabular}{|l|l|l|}
\multicolumn{3}{c}{\large Table 1}\\ [5pt]
\hline
$n$ &  $\lambda'$   &$\beta$\\ \hline
31& 3.06$K_{0}$ & 4.05$K_{0}$ \\ \hline
30-19& 4.6$K_{0}$ & 10$K_{0}$ \\ \hline
18&3.5$K_{0}$&  3.78$K_{0}$ \\ \hline
17& 3.26$K_{0}$&  3.60$K_{0}$ \\ \hline
16& 3.13$K_{0}$&  3.78$K_{0}$ \\ \hline
15& 2.76$K_{0}$&  3.12$K_{0}$ \\ \hline
14&2.34$K_{0}$&  2.96$K_{0}$ \\ \hline
13&2.03$K_{0}$&  2.15$K_{0}$ \\ \hline
\end {tabular}
$$

\end{document}